\documentclass[11pt]{article}

\usepackage{xypic,amsmath,amssymb,latexsym,theorem,enumerate,geometry,color,epic,eepic}
{\theorembodyfont{\rmfamily}\newtheorem{example}{Example}[section]}
{\theorembodyfont{\rmfamily}\newtheorem{rem}[example]{Remark}}
{\theorembodyfont{\rmfamily}}
{\theorembodyfont{\rmfamily}}
{\theorembodyfont{\rmfamily}}

\newcommand{\llabto}[2]{\stackrel{#2}
{\rule[0.54ex]{#1 em}{0.07ex}\hspace{-0.4em}\longrightarrow}}
\newcommand{\labto}[1]{\stackrel{#1}{\longrightarrow}}
\usepackage[all]{xy}

\def\leq{\leqslant}
\def\geq{\geqslant}
\newcommand{\xdirects}[2]{\def\objectstyle{\scriptstyle}
\objectmargin={0pt} \xy
(0,0)*+{}="a",(0,-6)*+{\rule{0em}{1.5ex}#2}="b",(7,0)*+{\;#1}="c"
\ar@{->} "a";"b" \ar @{->}"a";"c" \endxy}

\newcommand{\directs}[2]{\def\objectstyle{\scriptstyle}
\objectmargin={0pt} \xy
(0,0)*+{}="a",(0,-6)*+{\rule{0em}{1.5ex}#2}="b",(7,0)*+{\;#1}="c"
\ar@{->} "a";"b" \ar @{->}"a";"c" \endxy}
\newcommand{\tdirects}[2]{\def\objectstyle{\scriptstyle} \objectmargin={0pt}
\xy
(0,4)*+{}="a",(0,-2)*+{\rule{0em}{1.5ex}#2}="b",(7,4)*+{\;#1}="c"
\ar@{->} "a";"b" \ar @{->}"a";"c" \endxy }

\newcommand{\threeaxes}[3]{\def\objectstyle{\scriptstyle}  \objectmargin={0pt}
\xy
(0,0)*+{}="a",(0,-6)*+{\rule{0em}{1.5ex}#2}="b",(7,0)*+{\;#1}="c",
(14,-3)*+{\;#3}="d" \ar@{->} "a";"b" \ar @{->}"a";"c"  \ar
@{->}"a";"d"\endxy }

\newcommand{\threeaxesb}[3]{\def\objectstyle{\scriptstyle}  \objectmargin={0pt}
\xy
(0,0)*+{}="a",(0,-8)*+{\rule{0em}{1.5ex}#2}="b",(10,0)*+{\;#1}="c",
(8,7)*+{\;#3}="d" \ar@{->} "a";"b" \ar @{->}"a";"c"  \ar
@{->}"a";"d"\endxy }

\def\Rp{\mathbb{R}^+}
\def\epsilon{\varepsilon}
\def\del{\partial}
\def\A{\alpha}
\def\eps{\varepsilon}
\def\ge{\geq}
\def\le{\leq}

\begin{document}
\title{Moore hyperrectangles on a space form a  strict cubical omega-category}
\author{Ronald Brown}
\maketitle
\begin{abstract}
A question of Jack Morava is answered by generalising the notion of
Moore paths to that of Moore hyperrectangles, so obtaining a strict
cubical $\omega$-category. This also has the structure of
connections in the sense of Brown and Higgins, but cancellation of
connections does not hold.
\end{abstract}

\section*{Introduction}\label{sec:inro}
We recall in Section \ref{sec:paths} the notion of the space of
Moore paths on a topological space $X$.  A variant of the definition
is given in \cite{B2006}. Moore paths have the advantage of giving a
category of paths, with an associative composition, and identities,
rather than the common description in terms of maps $I \to X$.

However, whereas in higher dimensions the appropriate and analogous
operations on maps of cubes $I^n \to X$ have been well used, see for
example \cite{BH81:col}, there seems to have been in higher
dimensions no  definition analogous to that of Moore paths.

In this paper we give such a definition in Section \ref{sec:rect}
and in Section \ref{sec:laws} we give the laws that this structure
satisfies. The formulation of these is taken from \cite{ABS}, but
for a large part they go back to \cite{BH77}.

The cubical laws  were given in \cite{Kan55}. The cubical approach
in that paper was abandoned in favour of simplicial sets once the
problems of the geometric realisation of the cartesian product were
found, and Milnor had written on the geometric realisation of the
simplicial sets.

The introduction of `connections' in all dimensions was in
\cite{BH77,BH81:algcub,BH81:col} for the purpose of discussing
`commutative shells'. This was extended to the category case in
\cite{Mosa-phd,al-agl-phd,ABS}. The general theory of cubical sites
is developed in \cite{Grandis-cubsite}. Maltsiniotis has shown in
\cite{maltsin-cubandconn} that cubical sets with connections, in
contrast to the standard case,  have good realisations of cartesian
products. The thesis \cite{patchkoria-diplome} uses cubical sets
with what he calls {\it pseudo connections} for  the theory of
derived functors, analogously to the simplicial case.

The paper \cite{grandis-imageanal} uses an analogous procedure to
this for the definition of a higher categorical structure, with
cubes indexed on $Z^n$ and constant on each variable outside of a
certain `support', but does not take the support as part of the
structure.

An application of the cubical classifying space of a crossed complex
is in \cite{FennRouSan-trunks}.

It is suggested in recent work that the notion of Kan simplicial set
can be regarded as an $\infty$-groupoid, see for example
\cite{lurie-highertopos}. In some ways this is curious as this is
regarded as the start of such an idea is the fundamental group or
groupoid, which is made of classes of paths under homotopy  relative
to the end points. One would expect on the same principle to take
some form of homotopy classes of maps of $m$-paths. The difficulty
in this is shown that by the fact that an absolute strict homotopy
$m$-groupoid has been defined only for $m=1,2$, in \cite{BHKP}. The
paper \cite{BH81:col} shows that successful higher homotopy
groupoids can be defined for filtered spaces. This allows a route
into algebraic topology without setting up singular homology theory.

In any case, the construction $M_*(X)$ can be seen as another
candidate for a weak form of $\infty$-groupoid.

\section{Moore paths}\label{sec:paths}

Let  $\Rp=[0, \infty)$ be the nonnegative real line. For a space $X$
let $M(X)$ be the subspace of $X^{\Rp}\times \Rp$ of pairs $(f,r)$
such that $f$ is constant on $[r, \infty)$. There are two
maps\begin{align*}
\partial^-, \partial^+: M(X) &\to X,\\
\partial^-(f,r)&=f(0),\\
\partial^+(f,r)&=f(r).
\end{align*}
Now composition $\circ$ of Moore paths on $M(X)$ is given by the
composition
$$M(X)\, {}_{\partial^+}\hspace{-0.4em}\times _{\partial^-}M(X) \labto{\phi} X^{\Rp} \times \Rp \times \Rp \llabto{1}{1 \times + }X^{\Rp} \times \Rp  $$
where the first term is the pullback, and $\phi$   sends pairs
$(f,r),(g,s)\in M(X)$ such that $f(r)=g(0)$ to   triples $(h,r,s)\in
X^{\Rp}\times \Rp \times \Rp$ such that $h$ is constant on $[r+s,
\infty)$, $h|[0,r]=f|[0,r]$ and $h(t)=g(t-r)$ for $t \geq r$, and
$+$ is the addition function. So composition is continuous.

We also have an identity function $\epsilon: X \to M(X)$ given by
$\epsilon(x)=(\hat{x},0)$ where $\hat{x}$ is the constant map on
$\Rp$ with value $x$.

This composition gives, as is well known, a category structure
$(M(X), \partial^\pm,\circ,\epsilon)$. This structure also has a
`reverse' $-:M(X) \to M(X)$ given by $-(f,r)= (g,r)$ where $$g(t) =
\begin{cases}
  f(r-t) & \text{ if } 0 \leq t \leq r, \\
  f(0) &\text{ if } t \geq r.
\end{cases}$$
Thus $\partial^-(-a)=\partial^+(a), \partial^+(-a) = \partial^- a$.

We now discuss the relation with the fundamental groupoid on a set
$C$ of base points in $X$.

By a homotopy $H$ of elements $a^0=(f^0,r^0),a^1=(f^1,r^1)$ of
$M(X)$ we mean a continuous map $H: [0,1] \to M(X)$ such that
$H(0)=a^0, H(1)=a^1$, or, equivalently, a map $H: [0,1] \times \Rp
\to X $ such that $H(0,t)= a^0(t), H(1,t)=a^1(t)$ for $t \in \Rp$
and there is a continuous function $s \mapsto r(s)$ where $0 \leq s
\leq 1, r(s) \in \Rp$, $r(0)=r^0, r(1)=r^1$ and $H(s,t)= H(s,r(s))$
for $t \geq r(s), 0 \leq s \leq 1$. This homotopy is {\it rel end
points} if $H(s,0)=f^0(0), H(s,r(s))=f^0(r^0)$ for all $0 \leq s
\leq 1$. The fundamental groupoid $\pi_1(X,C)$ on the set of base
points $C \subseteq X$ is the the set of homotopy classes rel end
points of elements of $M(X)$ with source and target in $C$. For more
information on the use of $\pi_1(X,C)$, but with a slightly
different construction, see \cite{B2006}.

\section{Moore hyperrectangles}\label{sec:rect}

Let $M_n(X)$ be the subspace of $X^{(\Rp) ^n}\times (\Rp) ^n$ of
pairs $(f,(r))$ where $(r)=(r_1, \ldots, r_n)$ such that
$$f(t_1, \ldots, t_i, \ldots, t_n) = f(t_1, \ldots,r_i, \ldots, t_n) \text{ for } t_i \geq r_i, i=1, \ldots, n. $$
We call $(r)$ the {\it shape} and $f$ the {\it action}  of the
$n$-path $(f,(r))$.  We have
$$\partial^-_i,
\partial^+_i: M_n(X)\to M_{n-1}(X)$$ given by evaluating at $0$ or
$r_i$ in the $i$th position and omitting the $r_i$. More precisely,
$\partial^\alpha_i(f,(r))= (f',(r'))$ where $(r')=(r_1, \ldots,
\hat{r}_i, \ldots, r_n)$ and $f'(r')=f(r_1, \ldots, \alpha', \ldots,
r_n)$ where $\alpha'=0$ or $r_i$ according as $\alpha= -$ or $+$.

To define the degeneracies $\epsilon_i:M_{n-1}(X) \to M_{n}(X)$ we
set $\epsilon_i(f',(r'))=(f,(r))$ where $(r)$ is obtained from
$(r')$ by putting $0$ in the $i$th place, and $f(t_1, \ldots,
t_n)=f'(t_1, \ldots, \hat{t}_i, \ldots, t_n)$.

To define the connections $\Gamma^-_i:M_{n-1}(X) \to M_{n}(X)$ we
set $\Gamma^-_i(f',(r'))=(f,(r))$ where $(r)$ is obtained from
$(r')$ by repeating $r_i$ (in the $i$th and $(i+1)$th place, and
moving the others along), and setting $$f(t_1, \ldots, t_n) =f'(t_1,
\ldots,t_{i-1}, \max(t_i,t_{i+1}),t_{i+2}, \ldots, t_n).$$ Similarly
we get $\Gamma^+_i$ using $\min$ instead of $\max$. (This follows
the conventions of \cite{ABS}.)

For $i= 1, \ldots,n$ the category structure   $(M_n(X),
,\partial^-_i,\partial^+_i, \circ_i,\epsilon_i)$ is simply  that
given in section 1 but in the $i$th place.

In this way we give the family $M_*(X)=\{M_n(X)\}$ for $n \geq 0$
the structure of cubical  $\omega$-category: the laws for this and
the connections are given in Section \ref{sec:laws}. The paper
\cite{BH81:inf} also shows how to obtain what we now call a globular
$\omega$-category from this cubical structure, as a substructure in
which certain faces of a cube have various levels of degeneracy.
However this globular structure is not equivalent to the cubical
structure, as the proof in \cite{ABS} requires the cancellation law
for connections, which does not hold here: see Remark
\ref{rem:cancellationconn}.

We refer also to \cite{brownglobularhhgpd} for the construction of a
fundamental globular $\omega$-groupoid $\rho(X_*)$ of a filtered
space $X_*$.
\section{Laws} \label{sec:laws}
In this section we give the full structure and laws on the cubical
set with connections and compositions $M_*(X)$. We take these from
\cite{ABS}.

Let  $K$  be a cubical set, that is, a family of  sets  $\{ K_n;n
\ge 0\} $ with for $n \ge 1 $ face maps  $ \del _i^\A:K_n  \to
K_{n-1} \; (i = 1,2,\ldots,n;\, \alpha = +,-)$  and degeneracy maps
$\epsilon_i:K_{n-1}     \to   K_n \; (i = 1,2,\ldots,n)$ satisfying
the usual cubical relations:

\begin{alignat*}{2}
         \del_i^\A \del_j^\beta &=  \del_{j-1}^\beta \del_i^{\alpha}
        &&\hspace{-5em}(i<j),
        \tag*{(3.1)(i)} \\
        \epsilon_i \epsilon _j &=  \epsilon _{j+1} \epsilon _i && \hspace{-5em}(i \le j),
          \tag*{(3.1)(ii)} \\
\del ^{\alpha}_i \epsilon _j &=
                  {\begin{cases} \eps_{j-1} \del _i^\A &\hspace{8em} (i<j)  \\
                          \eps_{j} \del _{i-1}^\A & \hspace{8em}(i>j)  \\
                             \mathrm{id} &\hspace{8em} (i=j)
                  \end{cases}} && \tag*{(3.1)(iii)}   \\
\intertext{ We say that  $K$  is a  {\em cubical  set with
connections}  if  for $n \ge 0$ it  has additional structure maps
(called {\em connections}) $\Gamma_i^+,\Gamma_i^- :K_{n} \to K_{n+1}
\; (i = 1,2,\ldots,n)$ satisfying the relations:}
       \Gamma_i^\A   \Gamma_j^\beta & =  \Gamma_{j+1}^\beta\Gamma_i^\A &&
        \hspace{-5em}  (i  <  j)
       \tag*{(3.2)(i)} \\
         \Gamma_i^\A   \Gamma_i^\A  & =  \Gamma_{i+1}^\A\Gamma_i^\A &&
       \tag*{(3.2)(ii)} \\
       \Gamma_i^\A   \epsilon_j & = {\begin{cases}  \epsilon_{j+1}\Gamma_i^\A &\hspace{8em} (i < j)\\
                                     \epsilon_{j}\Gamma_{i-1}^\A &\hspace{8em}(i > j)
                                  \end{cases}}&& \tag*{(3.2)(iii)} \\
        \Gamma_j^\A \eps_j &= \eps^2_j=\eps_{j+1}\eps_j, &&  \tag*{(3.2)(iv)} \\
        \del^\A _i \Gamma_j^\beta &= {\begin{cases} \Gamma_{j-1}^\beta\del^\A_i
        & \hspace{8em}(i<j) \\
                              \Gamma_{j}^\beta\del^\A_{i-1}  &\hspace{8em} (i> j+1),
                              \end{cases}}&& \tag*{(3.2)(v) }\\
        \del^\A_j\Gamma_j^\A&= \del_{j+1}^\A   \Gamma _j^\A  =  id, && \tag*{(3.2)(vi)} \\
        \del^\A_j \Gamma^{-\A} _j&=
        \del^\A_{j+1}  \Gamma^{-\A} _j  = \epsilon_j  \del^\A_j.
        &&        \tag*{(3.2)(vii)}
\end{alignat*}
The connections are to be thought  of  as  extra  `degeneracies'. (A
degenerate cube of type  $ \epsilon_j  x$  has a pair of opposite
faces  equal and all other faces degenerate.  A cube of type  $
\Gamma_i^\A  x$  has a pair  of adjacent faces equal and all other
faces of type $\Gamma_j^\A  y$  or $\epsilon_j y$ .) Cubical
complexes with this, and other, structures  have  also been
considered by Evrard \cite{Ev}.

The prime example of  a  cubical   set  with  connections  is the
singular cubical complex  $KX$  of a space   $X$.   Here for $n \ge
0$ $K_n$ is  the  set  of singular $n$-cubes in $X$ (i.e. continuous
maps $I^n   \to   X$)  and the connection $ \Gamma_i^\A :K_{n } \to
K_{n+1}$  is induced by the map $\gamma_i^\A   : I^{n+1} \to
 I^{n}$    defined by
   $$   \gamma _i^\A (t_1 ,t_2 ,\ldots,t_{n+1} ) =
             (t_1 ,t_2 ,\ldots,t_{i-1},A(t_i ,t_{i+1}),t_{i+2},\ldots,t_{n+1} )
             $$
 where $A(s,t)=\max(s,t), \min(s,t)$ as $\A=-,+$ respectively.
 Here are pictures of $\gamma^\alpha_1 : I^2 \to I^1$ where the
 internal lines show lines of constancy of the map on $I^2$.


 \begin{center}
\setlength{\unitlength}{0.12in}
\begin{picture}(0,-10)(0,3)
\put(-5,2){\makebox{$\gamma^{-}_1 = $}}
\put(-17,2){\makebox{$\gamma^{+}_1 = $}}

\put(0,0){\framebox(4,4){}} \put(0,1){\line(2,0){3}}
\put(0,2){\line(3,0){2}} \put(2,4){\line(0,-2){2}}
\put(0,3){\line(3,0){1}} \put(2,4){\line(0,-2){2}}
\put(1,4){\line(0,-2){1}} \put(3,4){\line(0,-2){3 }}
\put(5,2){\tdirects{1}{2}}

\put(-12,0){\framebox(4,4){}} \put(-11,3){\line(2,0){3}}
\put(-11,3){\line(0,-1){3}} 
\put(-10,2){\line(3,0){2}} \put(-10,2){\line(0,-2){2}}
\put(-9,1){\line(0,-2){1}} \put(-9,1){\line(1,0){1 }}
\end{picture}
\end{center}
\vspace{1cm}

The complex  $KX$  has some further relevant structure, namely the
composition of $n$-cubes in the $n$  different directions.
Accordingly, we define a {\it cubical complex with connections and
compositions} to be a cubical set  $K$  with connections in which
each $K_n$ has  $n$ partial compositions $\circ _j\; (j =
1,2,\ldots,n)$ satisfying the following axioms. If  $a,b \in K_n$,
then  $a\circ _j b$  is defined if and only if $\del^-_j   b =
\del^+_j  a$  , and then
\begin{equation} \begin{cases} \del^-_j  (a\circ _j b) = \del^-_ja & \\
                 \del^+_j  (a\circ _j b) = \del^+_jb & \end{cases}
                 \qquad
 \del^\A_i  (a\circ _j b) =  \begin{cases} \del^\A_ja\circ _{j-1}\del^\A_i b &(i<j) \\
                 \del^\A_i  a\circ _j \del^\A_i b& (i>j), \end{cases}
                 \tag*{(3.3)}    \end{equation}

 {\em The  interchange laws}.  If  $i \ne j$  then
\begin{equation}
      (a\circ _i b) \circ _j  (c\circ _i d) = (a\circ _j c) \circ _i  (b\circ _j d)
      \tag*{(3.4)}
\end{equation}
whenever both sides are defined. (The diagram
$$
  \begin{bmatrix}
               a &b \\c&d
  \end{bmatrix}     \quad \tdirects{i}{j}
$$ will be used to indicate that both sides of the above equation
are  defined and also to denote the unique composite of the four
elements.)

     If  $i \ne j$  then
\begin{align*}
           \epsilon_i(a\circ _j b) &= \begin{cases}
           \eps_ia \circ _{j+1} \eps_ib & (i \le j) \\
           \eps_ia \circ _j\eps_ib & (i >j) \end{cases} \tag*{(3.5)} \\
       \Gamma^\A _i (a\circ _j b)& =  \begin{cases}
           \Gamma^\A_ia \circ _{j+1} \Gamma^\A_ib & (i < j) \\
           \Gamma^\A_ia \circ _j\Gamma^\A_ib & (i >j) \end{cases}
           \tag*{(3.6)(i)} \\
       \Gamma^+_j(a\circ _jb)&=  \begin{bmatrix}\Gamma^+_ja & \eps_j a\\
       \eps_{j+1} a & \Gamma^+_j b \end{bmatrix} \quad \tdirects{j}{j+1}  \tag*{(3.6)(ii)}\\
       \Gamma^-_j(a\circ _jb)&=  \begin{bmatrix}\Gamma^-_ja & \eps_{j+1} b\\
       \eps_{j} b & \Gamma^-_j b \end{bmatrix}\quad \tdirects{j}{j+1} \tag*{(3.6)(iii)} \\
\end{align*}        These last two equations are the {\it transport
laws}\footnote{Recall from \cite{BS76} that the term {\it
connection} was chosen because of an analogy with path-connections
in differential geometry. In particular, the transport law is a
variation or special case of the transport law for a
path-connection. }.

It is easily verified that the cubical Moore  complex   $M_*X$   of
a space $X$ satisfies these axioms with our above definitions. In
this context the transport law for $\Gamma^-_1 (a\circ b)$ can be
illustrated by the picture 
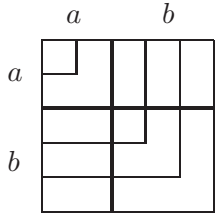
\vspace{0.9in}

\begin{center}
\begin{picture}(10,10)(40,20)

\setlength{\unitlength}{0.18in}
 \put(0,1){\line(4,0){5}} \put(0,1){\line(0,4){5}}
 \put(0,3){\line(4,0){3}}  \put(5,1){\line(0,4){5}}
 \put(0,2){\line(4,0){4}} \put(4,2){\line(0,4){4}}
\put(0,5){\line(4,0){1}}  \put(3,3){\line(0,4){3}}
\put(0,6){\line(4,0){5}}   \put(0,4){\line(4,0){5}}
\put(1,5){\line(0,4){1}}   \put(2,1){\line(0,4){5}}
\put(0,4.06){\line(4,0){5}} \put(2.06,1){\line(0,4){5}}
\put(0.7,6.5){\makebox{$a$}}
 \put(3.5,6.5){\makebox{$b$}}
 \put(-1,2.2){\makebox{$b$}}
\put(-1,4.8){\makebox{$a$}}
\end{picture}
\end{center}

\begin{rem}
That the above laws for  these structures apply to $M_*(X)$ is easy
to check. It is important that the shape tuples $(r_1, \ldots, r_n)$
are part of the structure. Thus if $\partial^+_1(f, (r))=
\partial^-_1(g,(s))$ then this implies that $r_i=s_i, 2 \leq i \leq
n$ as well as
$$f(r_1,t_2,\ldots,t_n)=g(0,t_2,\ldots,t_n) \text{ for } 0 \leq t_i \leq r_i, \; 2 \leq i \leq n. $$
This may seem a strong condition, but `composition is the inverse of
subdivision', and this enables one to obtain multiple compositions
as the inverse of `subdividing'  an element $(f,(r)) \in M_n(X)$.
\hfill $\Box$
\end{rem}
\begin{rem}\label{rem:cancellationconn}
In \cite{ABS} a cubical $\omega$-category with connections $G =
\{G_n\}$  is defined as a cubical set with connections and
compositions such that each  $\circ _j$   is a category structure on
$G_n$    with identity elements $\epsilon_j  y \;(y \in G_{n-1} )$,
and in addition
\begin{equation} \Gamma^+_ix \circ_i\Gamma^-_ix = \eps_{i+1}x, \quad
\Gamma^+_ix \circ_{i+1}\Gamma^-_ix = \eps_{i}x.\tag{2.7}
\end{equation}
However this cancellation law does not hold for $M_*(X)$. Thus the
equivalence between globular and cubical categories developed in
\cite{ABS} does not apply to $M_*(X)$, nor does  the exact relations
between `commutative shells' and `thin elements' developed in
\cite{higgins-thin}. \rule{1cm}{0ex}\hfill $\Box$
\end{rem}
\begin{rem}
There are  also {\it reverses} $-_i: M_n(X) \to M_n(X), 1 \leq i
\leq n$ defined as in Section \ref{sec:paths}. A problem with our
construction is that a path $[0,1] \to M_n(X)$ is not necessarily an
element  of $M_{n+1}(X)$. In particular the easily defined homotopy
rel end points of paths $a \circ -a \simeq 0_{\partial^-a}$ is not
an element of $M_2(X)$. \hfill $\Box$
\end{rem}
\section{Tensor products}
\label{sed:tensor} The tensor product $K \otimes L$ of cubical sets
is also defined in \cite{Kan55} and extended to cubical sets with
connections in \cite{BH87,ABS}. We see the convenience of cubes in
the current account since since if $a=(f,(r)) \in M_m(X)$ and
$b=(g,(s)) \in M_n(Y)$ then their tensor product $a \otimes b\in
M_{m+n}(X \times Y)$. It is given by $a \otimes b=(h,((r)\circ(s)))$
where $(r) \circ (s)= (r_1, \ldots, r_m, s_1, \ldots,s_n)$ and $h=f
\times g : (\Rp)^m \times (\Rp)^n \to X \times Y$ with the usual
identification of $(\Rp)^m \times (\Rp)^n $ and $(\Rp)^{m+n}$.

\end{document}